\documentclass[11pt]{amsart}
\usepackage{amsmath,amsfonts}

\newtheorem{thm}{Theorem}[section]
\newtheorem{cor}[thm]{Corollary}
\newtheorem{lem}[thm]{Lemma}

\newtheorem{rem}[thm]{Remark}

\newtheorem{op}[thm]{Open Problem}

\numberwithin{equation}{section}

\begin{document}

\title[Rolle's theorem is either false or trivial]
{Rolle's theorem is either false or trivial in
infinite-dimensional Banach spaces}
\author{Daniel Azagra and Mar Jim\'enez-Sevilla}

\dedicatory{Dedicated to Albert Galvany and Pilar Olivella, who among other
residents of the Colegio de Espa\~na in Paris contributed to
create the mood of magic rationality in which the tubes were
conceived.}


\begin{abstract}
We prove the following new characterization of $C^p$ (Lipschitz)
smoothness in Banach spaces. An infinite-dimensional Banach space
$X$ has a $C^p$ smooth (Lipschitz) bump function if and only if it
has another $C^p$ smooth (Lipschitz) bump function $f$ such that
$f'(x)\neq 0$ for every point $x$ in the interior of the support
of $f$ (that is, $f$ does not satisfy Rolle's theorem). Moreover,
the support of this bump can be assumed to be a smooth starlike
body. As a by-product of the proof of this result we also obtain
other useful characterizations of $C^p$ smoothness related to the
existence of a certain kind of deleting diffeomorphisms, as well
as to the failure of Brouwer's fixed point theorem even for smooth
self-mappings of starlike bodies in all infinite-dimensional
spaces. Finally, we study the structure of the set of gradients of
bump functions in the Hilbert space $\ell_2$, and as a consequence
of the failure of Rolle's theorem in infinite dimensions we get
the following result. The usual norm of the Hilbert space $\ell_2$
can be uniformly approximated by $C^1$ smooth Lipschiz functions
$\psi$ so that the cones generated by the sets of derivatives
$\psi'(\ell_{2})$ have empty interior. This implies that there are
$C^1$ smooth Lipschitz bumps in $\ell_{2}$ so that the cones
generated by their sets of gradients have empty interior.
\end{abstract}

\maketitle

\section[Introduction]{Introduction and main results}

Rolle's theorem in finite-dimensional spaces states that, for
every bounded open subset $U$ of $\mathbb{R}^{n}$ and for every
continuous function $f:\overline{U}\longrightarrow\mathbb{R}$
such that $f$ is differentiable in $U$ and constant on the
boundary $\partial U$, there exists a point $x\in U$ such that
$f'(x)=0$. Unfortunately, Rolle's theorem does not remain valid
in infinite dimensions. It was S. A. Shkarin \cite{Shk} that
first showed the failure of Rolle's theorem in superreflexive
infinite-dimensional spaces and in non-reflexive spaces which
have smooth norms. The class of spaces for which Rolle's theorem
fails was substantially enlarged in \cite{AGJ}, where it was also
shown that an approximate version of Rolle's theorem remains
nevertheless true in all Banach spaces. In fact, as a consequence
of the existence of diffeomorphisms deleting points in
infinite-dimensional spaces (see \cite{A, Ado}), it is easy to
see that Rolle's theorem fails in all infinite-dimensional Banach
spaces which have smooth norms \cite{thesis}.

However, none of these results allows to characterize the spaces
for which Rolle's theorem fails since, as shown by R. Haydon
\cite{Haydon}, there are Banach spaces with smooth bump functions
which do not possess any equivalent smooth norms. Of course,
Rolle's theorem is trivially true in the Banach spaces which do
not have any smooth bumps (if $X$ is such a space then every
function on $X$ satisfying the hypothesis of Rolle's theorem must
be a constant). Thus, in many infinite-dimensional Banach spaces,
Rolle's theorem either fails or is trivial, depending on the
smoothness properties of the spaces considered. In this setting,
it does not seem too risky to conjecture, as it was done in
\cite{AGJ}, that Rolle's theorem should fail in an
infinite-dimensional Banach space if and only if our space has a
$C^1$ smooth bump function. In this paper we will prove this
conjecture to be right, thus providing an interesting new
characterization of smoothness in Banach spaces.

Our main result is the following

\begin{thm}\label{Rolle's theorem fails}
Let $X$ be an infinite-dimensional Banach space which has a $C^p$
smooth (Lipschitz) bump function. Then there exists another $C^p$
smooth (Lipschitz) bump function $f:X\longrightarrow [0,1]$ with
the property that $f'(x)\neq 0$ for every $x\in\textrm{int}(\textrm{supp} f)$.
\end{thm}
\noindent Here $\textrm{supp} f$ denotes the support of $f$, that
is, $\textrm{supp} f=\overline{\{x\in X : f(x)\neq 0\}}$. Let us
recall that $b:X\longrightarrow\mathbb{R}$ is said to be a bump
function on $X$ provided $b$ is not constantly zero and $b$ has a
bounded support.

\vspace{0.2cm} From this result it is easily deduced the following

\begin{cor}\label{Rolle's failure characterizes smoothness}
Let $X$ be an infinite-dimensional Banach space. The following
statements are equivalent.
\begin{itemize}
\item [{(1)}] $X$ has a $C^p$ smooth (and Lipschitz) bump function.
\item [{(2)}] There exist a bounded contractible open subset $U$ of $X$ and a
continuous function $f:\overline{U}\longrightarrow\mathbb{R}$
such that $f$ is $C^p$ smooth (and Lipschitz) in $U$,  $f=0$ on
$\partial U$, and yet $f'(x)\neq 0$ for all $x\in U$, that is,
Rolle's theorem fails in $X$.
\item [{(3)}] There exist a $C^p$ smooth (and Lipschitz) function
$f:X\longrightarrow [0,1]$ and a bounded contractible open subset
$U$ of $X$ such that $f=0$ precisely on $X\setminus U$ and yet
$f'(x)\neq 0$ for all $x\in U$.
\end{itemize}
\end{cor}

Just in order to complete the picture of Rolle's theorem in
infinite-dimensional Banach spaces let us quote the two positive
results from \cite{AGJ, Adev1} on approximate and subdifferential
substitutes of Rolle's theorem, which guarantee the existence of
arbitrarily small derivatives (instead of vanishing ones) for
every function sayisfying (in an approximate manner) the
conditions of the classic Rolle's theorem.

\begin{thm}[Azagra--G\'omez--Jaramillo]\label{Approximate Rolle's
theorem} Let $U$ be a bounded connected open subset of a Banach
space $X$. Let $f:\overline{U}\longrightarrow\mathbb{R}$ be a
bounded continuous function which is (G\^ateaux) differentiable in
$U$. Let $R>0$ and $x_{0}\in U$ be such that $\textrm{dist}(x_{0},
\partial U)=R$. Suppose that $f(\partial U)\subseteq [-\varepsilon, \varepsilon]$
for some $\varepsilon>0$. Then there exists some
$x_{\varepsilon}\in U$ such that
$\|f'(x_{\varepsilon})\|\leq\frac{\varepsilon}{R}$.
\end{thm}

\begin{thm}[Azagra--Deville]\label{Approximate Subdifferential Rolle's
theorem} Let $U$ be a bounded connected open subset of a Banach
space $X$ which has a $C^1$ smooth Lipschitz bump function. Let
$f:\overline{U}\longrightarrow\mathbb{R}$ be a bounded continuous
function, and let $R>0$ and $x_{0}\in U$ be such that
$\textrm{dist}(x_{0},
\partial U)=R$. Suppose that $f(\partial U)\subseteq [-\varepsilon, \varepsilon]$
for some $\varepsilon>0$. Then
    $$
    \inf\{\|p\| : p\in D^{-}f(x)\cup D^{+}f(x), x\in
    U\}\leq\frac{2\varepsilon}{R}.
    $$
\em{(Here $D^{-}f(x)$ and $D^{+}f(x)$ denote the subdifferential
and superdifferential sets of $f$ at $x$, respectively; see
\cite{DGZ}, p. 339, for the definitions).}
\end{thm}

The ``twisted tube'' method that we develop in section 2 in order to prove theorem
\ref{Rolle's theorem fails} is interesting in itself and, with
little more work, provides a useful
characterization of $C^p$ smoothness in infinite-dimensional
Banach spaces related to the existence of a certain kind of
{\em deleting diffeomorphisms}. Namely, we have the following

\begin{thm}\label{Characterization by diffeomorphisms deleting twisted tubes}
Let $X$ be an infinite-dimensional Banach space. The following
assertions are equivalent.
\begin{itemize}
\item [{(1)}] $X$ has a $C^p$ smooth bump function.
\item [{(2)}] There exists a nonempty contractible closed
subset $D$ of the unit ball $B_X$
and a $C^p$ diffeomorphism $f:X\longrightarrow
X\setminus D$ so that $f$ restricts to the identity outside $B_X$.
\end{itemize}
\end{thm}
This result yields the following corollaries.

First, the celebrated Brouwer's fixed point theorem
fails even for smooth self-mappings of balls or
starlike bodies in all infinite-dimensional Banach spaces.
Let us recall that Brouwer's theorem states that every
continuous self-map of the unit ball of a finite-dimensional
normed space admits a fixed point. This is the same as saying
that there is no continuous retraction from the unit ball onto
the unit sphere, or that the unit sphere is not contractible (the
identity map on the sphere is not homotopic to a constant map).
In infinite dimensions the situation is completely different
and Brouwer's theorem is no longer true (see \cite{BP, Nowak, BS, LS, GK, BL,
AC}. Theorem \ref{Characterization by diffeomorphisms deleting twisted tubes}
yields a trivial proof that Brouwer's theorem is false in infinite
dimensions even for smooth self-mappings of balls or starlike bodies; this is
a particular case (the non-Lipschitz one) of the main result in
\cite{AC}.

Second, we deduce from the above characterization that the support of the bump
functions which violate Rolle's theorem can
always be assumed to be a smooth starlike body. This is all shown
in section 3.

In section 2 we  give the proofs of theorems \ref{Rolle's theorem
fails} and \ref{Characterization by diffeomorphisms deleting
twisted tubes}. A much simpler proof of theorem  \ref{Rolle's theorem
fails} for the non-Lipschitz case is included in this section too.

Finally, in section 4 we study the structure of the set of
gradients of bump functions in the Hilbert space $\ell_2$, and as
a consequence of the failure of Rolle's theorem in infinite
dimensions we get the following result. The usual norm of the
Hilbert space $\ell_2$ can be uniformly approximated by $C^1$
smooth Lipschitz functions $\psi$ so that the cones generated by
the sets of derivatives $\psi'(\ell_{2})$ have empty interior.
This implies that there are $C^1$ smooth Lipschitz bumps in
$\ell_{2}$ so that the cones generated by their sets of gradients
have empty interior.


\section[Proofs]{The proofs}

The idea behind the proof of theorem \ref{Rolle's theorem fails}
is as simple as this. First we build a twisted tube $T$ of infinite
length in the interior of the unit ball $B_X$, with a beginning
but with no end. This twisted tube can be thought of as directed
by an ever-winding infinite path $p$ that gets lost in the
infinitely many dimensions of our space $X$. In technical words,
one can construct a diffeomorphism $\pi$ between a straight
(unbounded) half-cilynder $C$ and a twisted (bounded) tube $T$
contained in $B_X$. The tube $T$ is going to be the support of a
smooth bump function $f$ that does not satisfy Rolle's
theorem. In order to define such a function $f$ we only have to
make it strictly increase in the direction which is tangent to the
leading path $p$ at each point of the tube $T$. The graph of $f$
would thus represent an ever-ascending stairway built upon our
twisted tube, with a beginning but no end.

The spirit of the proof that (1) implies (2) in theorem
\ref{Characterization by diffeomorphisms deleting twisted tubes}
is not very different.
We will make use of the diffeomorphism $\pi$ between a straight
(unbounded) half-cilynder $C$ and a bounded twisted tube $T$
contained in $B_X$. If we consider a straight closed
half-cilinder $C'$ contained in the interior of $C$ and directed
by the same line as $C$, it is elementary that there is a
diffeomorphism $g:X\longrightarrow X\setminus C'$ so that $g$
restricts to the identity outside $C$. In fact this is true even in the
plane.
Now, by composing this diffeomorphism $g$ with the
diffeomorphisms $\pi$ and $\pi^{-1}$ that give us an appropriate
coordinate system in the twisted tube $T=\pi(C)$, we get a
diffeomorphism $f:X\longrightarrow X\setminus T'$, where
$T'=\pi(C')$ is a smaller closed twisted tube inside $T$, and $f$
restricts to the identity outside the unit ball. The precise
definition of $f$ would be $f(x)=\pi(g(\pi^{-1}(x)))$ if $x\in
T$, and $f(x)=x$ if $x\in X\setminus T$. If we take $D=T'$ we are
done.

In the rest of his section we will be involved in the task of
formalizing these ideas.

\vspace{0.2cm}

The following theorem guarantees the existence of bounded infinite
twisted tubes in all infinite-dimensional Banach spaces.

\begin{thm}\label{Existence of twisted tubes}
There are universal constants $M>0$ (large) and $\varepsilon>0$
(small) such that, for every infinite-dimensional Banach space
$X$, if we consider the decomposition $X=H\oplus [z]$ (where
$H=\textrm{Ker}\, z^{*}$ for some $z^{*}\in X^*$ with
$z^{*}(z)=\|z^*\|=\|z\|=1$) and the open half-cilynder $C$ of
diameter $\varepsilon$, directed by $z$, and with base on $H$,
$C=\{x+tz\in X : \|x\|<\varepsilon, t>0\}$, then there exists an
injection $\pi: C\longrightarrow B_X$ which is a $C^{\infty}$
diffeomorphism onto its image. The image $T=\pi(C)$ is thus a
bounded open set which we will call a bounded open infinitely
twisted tube in $X$. Moreover, the derivatives of the mappings
$\pi:C\longrightarrow T$ and  $\pi^{-1}:T\longrightarrow C$ are
both uniformly bounded by $M$.
\end{thm}

Assume for a while  that theorem \ref{Existence of twisted tubes}
is already established and let us explain how theorems \ref{Rolle's
theorem fails} and \ref{Characterization by diffeomorphisms deleting
twisted tubes} can be deduced.

\begin{center}
{\bf Proof of theorem \ref{Rolle's theorem fails}}.
\end{center}

Consider the diffeomorphism $\pi:C\longrightarrow T\subset B_X$
from theorem \ref{Existence of twisted tubes}. Take a $C^p$
smooth (Lipschitz) non-negative bump function $\varphi$ on $H$ so
that the support of $\varphi$ is contained in the base of $C$,
that is, $\varphi(x)=0$ whenever
$\|x\|\geq\frac{\varepsilon}{2}$, for instance. Pick a
$C^{\infty}$ smooth real function $\mu:\mathbb{R}\longrightarrow
[0,1]$ such that $\mu(t)=0$ for $t\leq 1$, $0<\mu(t)<1$ for $t>1$
and $0<\mu'(t)<1$ for all $t>1$. Then define $g:X=H\oplus
[z]\longrightarrow\mathbb{R}$ by
    $$
    g(x,t)=\varphi(x)\mu(t).
    $$
It is plain that $g$ is a $C^p$ smooth (Lipschitz) function such
that $g'(x,t)\neq 0$ for every $x\in\textrm{int}(\textrm{supp} f)$, that
is, for every $x$ such that $g(x,t)\neq 0$ (take into account that the
interior of the support of $g$ coincides in this case with the open support
of $g$, that is the set of points at which $g$ does not vanish). Indeed,
    $$
    g'(x,t)(0,1)=\frac{\partial g}{\partial t}(x,t)=\varphi(x)\mu'(t)
    $$
and therefore $g'(x,t)(0,1)=0$ if and only if $\varphi(x)=0$ or
$\mu'(t)=0$, which happens if and only if $\varphi(x)=0$ or
$\mu(t)=0$, that is to say, $g(x,t)=0$. Now let us define
$f:X\longrightarrow\mathbb{R}$ by
    $$
    f(y)=\left\{
    \begin{array}{ll}
    &g(\pi^{-1}(y)) \quad \textrm{if $y\in T$;}\\
    &0 \quad \textrm{if $y\notin T$}
    \end{array}\right.
    $$
It is clear that $f$ is a well defined $C^p$ smooth (Lipschitz)
function, and $\textrm{supp}(f)=\pi (\textrm{supp}(g))\subset T$,
from which it follows that $f$ has a bounded support. We claim
that $f'(y)\neq 0$ whenever $y\in\textrm{int}(\textrm{supp} f)$,
that is, $f$ does not satisfy Rolle's theorem. Indeed,  if
$y\in\textrm{int}(\textrm{supp} f)$ then
$\pi^{-1}(y)=(x,t)\in\textrm{int}(\textrm{supp} g)$ and therefore
$g'(x,t)(0,1)\neq 0$. But then
    $$
    f'(y)=g'(x,t)\circ D\pi^{-1}(y)\neq 0,
    $$
because $D\pi^{-1}(y)$ is a linear isomorphism. This concludes the
proof of theorem \ref{Rolle's theorem fails}.


\vspace{0.3cm}

Now we will turn our attention to the proof of theorem
\ref{Characterization by diffeomorphisms deleting twisted tubes}.
Before proceeding with the proof, let us fix some standard terminology and notation used
throughout this section and the following one. A closed subset $A$ of
a Banach space $X$ is said to be a starlike body provided $A$ has
a non-empty interior and there exists a point $x_{0}\in
\textrm{int}A$ such that each ray emanating from $x_{0}$ meets the
boundary of $A$ at most once. In this case we will say that $A$ is
{\em starlike with respect to $x_{0}$}. When dealing with starlike
bodies, we can always assume that they are starlike with respect
to the origin (up to a suitable translation), and we will do so
unless otherwise stated.

For a starlike body $A$, the characteristic cone of $A$ is defined
as $$ cc A=\{x\in X\mid r x\in A\hspace{2mm}\textrm{for all}
\hspace{2mm}r>0\}, $$ and the Minkowski functional of $A$ as $$
q_{A}(x)=\inf\{\lambda>0\mid \frac{1}{\lambda}x\in A\} $$ for all
$x\in X$. It is easily seen that for every starlike body $A$ its
Minkowski functional $q_{A}$ is a continuous function which
satisfies $q_{A}(rx)=rq_{A}(x)$ for every $r\geq 0$ and
$q_{A}^{-1}(0)=cc A$. Moreover, $A=\{x\in X \mid q_{A}(x)\leq
1\}$, and $\partial A=\{x\in X \mid q_{A}(x)=1\}$, where
$\partial A$ stands for the boundary of $A$. Conversely, if
$\psi:X\to [0, \infty)$ is continuous and satisfies $\psi(\lambda
x)=\lambda\psi(x)$ for all $\lambda\geq 0$, then $A_{\psi}=\{x\in
X \mid \psi(x)\leq 1\}$ is a starlike body. Convex bodies (that
is, closed convex sets with nonempty interior) are an important
kind of starlike bodies. We will say that $A$ is a $C^{p}$ smooth
(Lipschitz) starlike body provided its Minkowski functional
$q_{A}$ is $C^{p}$ smooth (and Lipschitz) on the set $X\setminus
q_{A}^{-1}(0)$.

It is worth noting that for every Banach
space $(X, \|.\| )$ with a $C^{p}$ smooth (Lipschitz) bump function there
exist a functional $\psi$ and constants $a, b>0$ such that $\psi$
is $C^{p}$ smooth (Lipschitz) away from the origin, $\psi(\lambda
x)=|\lambda|\psi(x)$ for every $x\in X$ and
$\lambda\in\mathbb{R}$, and $a\|x\|\leq \psi(x)\leq b\|x\|$ for
every $x\in X$ (see \cite{DGZ}, proposition II.5.1).
The level sets of this function are precisely
the boundaries of the smooth bounded starlike bodies $A_{c}=\{x\in
X \mid \psi(x)\leq c\}$, $c\in\mathbb{R}$.
This shows in particular that every Banach space having a $C^p$ smooth
(Lipschitz) bump function has a $C^p$ smooth (Lipschitz) bounded
starlike body as well. The converse is obviously true too.

\begin{center}
{\bf Proof of theorem \ref{Characterization by diffeomorphisms deleting twisted tubes}}.
\end{center}

First of all let us choose a
number $\varepsilon>0$, a cilynder $C$, a bounded twisted tube
$T$, and a diffeomorphism $\pi:C\longrightarrow T$
from theorem \ref{Existence of twisted tubes}.

Let $B$ be a $C^{\infty}$ smooth convex body in the plane
$\mathbb{R}^{2}$ whose boundary contains the set
    $$
    \{(s,t) : t=-1, |s|\leq\frac{\varepsilon}{4}\}
    \cup\{(s,t) : |s|=\frac{\varepsilon}{2},
    t\geq -1+\frac{\varepsilon}{4}\},
    $$
and let $q_B$ be the Minkowski functional of $B$. Define
$B'=\frac{1}{2}B=\{(s,t) : q_{B}(s,t)\leq\frac{1}{2}\}$. Let
$\theta:(\frac{1}{2},\infty)\longrightarrow [0,\infty)$ be a
$C^{\infty}$ smooth real function so that $\theta'(t)<0$ for
$\frac{1}{2}<t<1$, $\theta(t)=0$ for $t\geq 1$, and
$\lim_{t\to\frac{1}{2}^{+}}\theta(t)=+\infty$. Now define
$\varphi:\mathbb{R}^{2}\setminus B'\longrightarrow\mathbb{R}^{2}$
by
    $$
    \varphi(s,t)=(\varphi_{1}(s,t),\varphi_{2}(s,t))=
    (s, t+\theta(q_{B}(s,t))).
    $$
It is elementary to check that $\varphi$ is a $C^{\infty}$
diffeomorphism from $\mathbb{R}^{2}\setminus B'$ onto
$\mathbb{R}^{2}$ so that $\varphi$ restricts to the identity
outside the band $B$.

Next, recall that since $X$ has a $C^p$ smooth bump then it has a $C^p$
bounded starlike body $A$ as well.
If $X=H\oplus[z]$, take $W=A\cap H$, which is a $C^p$
bounded starlike body in $H$, and denote by $q_W$ its
Minkowski functional. We can assume that $W\subseteq B(0,1)$, that
is, $\|x\|\leq q_{W}(x)$ for all $x\in H$. Let us define
    $$
    \psi(x,t)=q_{B}(q_{W}(x),t)
    $$
for all $(x,t)\in X$. It is clear that $\psi$ is a continuous
function which is positive-homogeneous and $C^p$ smooth away from
the half-line $L=\{(x,t)\in X : x=0, t\geq 0\}$. Then the sets
    $$
    U=\{(x,t)\in X : \psi(x,t)\leq 1\}, \quad
    U'=\{(x,t)\in X : \psi(x,t)\leq\frac{1}{2}\}
    $$
are cilyndrical $C^p$ starlike bodies whose
characteristic cones are the half-line $L$. If we define
    $$
    h(x,t)=(x, \varphi_{2}^{-1}(q_{W}(x), t))
    $$
for $(x,t)\in X$, it is not difficult to realize that $h$ is a
$C^p$ diffeomorphism from $X$ onto $X\setminus U'$ so
that $h$ restricts to the identity outside $U$. The inverse of $h$
is given by
    $$
    h^{-1}(x,t)=(x, t+\theta(\psi(x,t))).
    $$

Now consider
the cilyndrical bodies $V:=(0,2)+U$ and $V':=(0,2)+U'$, and put
$g(x,t)=h(x, t-2)$. Then $g:X\longrightarrow X\setminus V'$ is a
$C^p$ diffeomorphism such that $g$ is the identity
outside $V$. Note that, since $W\subseteq B(0,1)$, we have that
$V'\subset V\subset C=\{(x,t)\in X : \|x\|<\varepsilon, t>0\}$.
Let us define
    $$
    f(x)=\left\{
    \begin{array}{ll}
    &\pi(g(\pi^{-1}(x))) \quad \textrm{if $x\in T$;}\\
    &x \quad \textrm{otherwise.}
    \end{array}\right.
    $$
It is then clear that $f$ is a $C^p$ diffeomorphism from
$X$ onto $X\setminus T'$, where $T'=\pi(V')$ is a smaller closed twisted
tube inside $\pi(V)\subseteq T$, and $f$ restricts to the identity
outside the larger tube $\pi(V)\subset T$, which is contained in $B_X$.
This completes the proof that (1) implies (2).

\medskip

Conversely, if there is such an $f$ as in (2), we can assume that $f(0)\neq 0$
and take $T\in X^*$ so that $T(f(0))\neq 0$; then the function
$b:X\longrightarrow\mathbb{R}$ defined by $b(x)=T(x-f(x))$ is a
$C^p$ smooth  bump on $X$.

\vspace{0.3cm}

Now we proceed with the  proof of theorem \ref{Existence of
twisted tubes}. We will make use of the following lemma, which
guarantees the existence of an appropiate path of linear
isomorphisms. Here $\textrm{Isom}(X)$ stands for  the set of
linear isomorphisms of $X$, which is regarded as a subset of
$\mathcal{L}(X,X)$, the linear continuous mappings of $X$ into
$X$.

\begin{lem}\label{Path beta of linear isomorphisms}
There is a universal constant $K>0$ such that for every
infinite-dimensional Banach space $X$ there are paths
$\beta:[0,\infty)\longrightarrow\textrm{Isom}(X)$ and
$p:[0,\infty)\longrightarrow X$ with the following properties:
\begin{itemize}
\item [{(1)}] Both $\beta$ and $p$ are $C^{\infty}$ smooth, as
well as the path of inverse isomorphisms
$\beta^{-1}:[0,\infty)\longrightarrow\textrm{Isom}(X)$,
$\beta^{-1}(t)=[\beta(t)]^{-1}$.
\item [{(2)}] $1\leq\|\beta(t)\|\leq K$ and $1\leq\|\beta^{-1}(t)\|\leq
K$ for all $t\in [0, \infty)$.
\item [{(3)}] $\sup_{t\geq 0}\|\beta'(t)\|\leq K$ and
$\sup_{t\geq 0}\|(\beta^{-1})'(t)\|\leq K$.
\item [{(4)}] There exists a certain $v\in X$, with
$1\geq\|v\|\geq\frac{1}{K}$, such that $p'(t)=\beta(t)(v)$ for all
$t\geq 0$.
\item [{(5)}] For every $t,s\in [0,\infty)$ we have that
$\|p(t)-p(s)\|\geq\frac{1}{K}\min\{1, |t-s|\}$.
\end{itemize}
\end{lem}
\begin{proof}

Let $(x_{n})_{n=0}^{\infty}$ be a normalized basic sequence in $X$
with biorthogonal functionals $(x_{n}^{*})_{n=0}^{\infty}\subset
X^*$ (that is, $x_{n}^{*}(x_{k})=\delta_{n,k}=1$ if $n=k$, and $0$
otherwise) satisfying $\|x_{n}^{*}\|\leq 3$ (one can always take
such sequences, see \cite{Day}, p. 93, or \cite{Diestel}, p. 39).
For $n\geq 1$ set
$v_{n}=x_{n}-x_{n-1}$. Let
$\theta:\mathbb{R}\longrightarrow\mathbb{R}$ be a $C^{\infty}$
function with the following properties:
\begin{itemize}
\item [{(i)}] $\theta(t)=0$ whenever $t\leq -\frac{1}{2}$ or
$t\geq 1$;
\item [{(ii)}] $\theta(t)=1$ for $t\in [0, \frac{1}{2}]$;
\item [{(iii)}] $\theta'(t)>0$ for $t\in (-\frac{1}{2},
0)$;
\item [{(iv)}] $\theta(t)=1-\theta(t-1)$ for $t\in [\frac{1}{2},
1]$;
\item [{(v)}] $\sup_{t\in\mathbb{R}}|\theta'(t)|\leq 4$.
\end{itemize}

For $n\geq 1$ let us define
$\theta_{n}:\mathbb{R}\longrightarrow\mathbb{R}$ by
$\theta_{n}(t)=\theta(t-n+1)$. It is clear that the functions
$\theta_n$ are all $C^{\infty}$ smooth and have Lipschitz constant
less than or equal to $4$, $\theta_{n}=0$ on $(-\infty,
n-1-\frac{1}{2}]\cup[n,\infty)$, $\theta_{n}=1$ on
$[n-1, n-\frac{1}{2}]$, and
$\theta_{n}(t)=1-\theta_{n+1}(t)$ for all $t\in [n-1,
n+\frac{1}{2}]$.

Our path $\beta$ of linear isomorphisms is going to be of the form

    $$
    \beta(t)=\sum_{n=1}^{\infty}\theta_{n}(t)S_{n},
    $$
where each $S_{n}\in\textrm{Isom}(X)$ takes the vector $v_{1}$
into $v_{n}$ and for every $\lambda\in[0, 1]$ the mapping
$L_{n,\lambda}=(1-\lambda)S_{n}+\lambda S_{n+1}$ is still a linear
isomorphism and, moreover, the families of isomorphisms
$\{L_{n,\lambda}\}_{n\in\mathbb{N}, \lambda\in [0,1]}$ and
$\{L^{-1}_{n,\lambda}\}_{n\in\mathbb{N}, \lambda\in [0,1]}$ are
uniformly bounded. Let us define the isomorphisms $S_n$. They are
going to be of the form
    $$
    S_{n}(x)=x+f_{n}(x)(v_{n}-v_{1}),
    $$
where $f_{n}\in X^*$ satisfies $f_{n}(v_{1})=1=f_{n}(v_{n})$, and
$\|f_{n}\|\leq 18$ (the exact definition of $f_n$ will be given
later). Their inverses $S_{n}^{-1}$ will be
    $$
    S_{n}^{-1}(y)=y-f_{n}(y)(v_{n}-v_{1}).
    $$
We want the linear mappings
$L_{n,\lambda}=(1-\lambda)S_{n}+\lambda S_{n+1}$ to be linear
isomorphisms. We have
    $$
    y=L_{n,\lambda}(x)=x+(1-\lambda)f_{n}(x)(v_{n}-v_{1})+
    \lambda f_{n+1}(x)(v_{n+1}-v_{1}),  \eqno (1)
    $$
from which
    $$
    x=y-[(1-\lambda)f_{n}(x)(v_{n}-v_{1})+
    \lambda f_{n+1}(x)(v_{n+1}-v_{1})],  \eqno (2)$$
and we need to write $f_{n}(x)$ and $f_{n+1}(x)$ as linear
functions of $y$. If we apply the functionals $f_n$ and $f_{n+1}$
successively to equation (1), we denote $A_{n}=f_{n}(x)$,
$B_{n}=f_{n+1}(x)$, $C_{n}=f_{n}(y)$, $D_{n}=f_{n+1}(y)$, and we
take into account that
$1=f_{n}(v_{1})=f_{n}(v_{n})=f_{n+1}(v_{1})$, then we obtain the
system
    $$
    \left\{
    \begin{array}{ll}
    &A_{n}+\lambda [f_{n}(v_{n+1})-1]B_{n}=C_{n}\\
    &(1-\lambda)[f_{n+1}(v_{n})-1]A_{n}+B_{n}=D_{n},
    \end{array}\right.   \eqno (3)
    $$
which we want to have a unique solution for $A_{n}, B_{n}$. The
determinant of this system is
    $$
    \Delta_{n,\lambda}=1-\lambda
    (1-\lambda)[f_{n+1}(v_{n})-1][f_{n}(v_{n+1})-1],
    $$
and we want $\Delta_{n,\lambda}$ to be bounded below by a strictly
positive number, and this bound has to be uniform in $n,\lambda$.
For $n\geq 3$ this can easily be done by setting
    $$
    f_{n}=x_{1}^{*}-x_{n-1}^{*}
    $$
(so that $f_{n}(v_{n})=1=f_{n}(v_{1})$, $f_{n}(v_{n+1})=0$,
$f_{n+1}(v_{n})=-1$, and therefore
$\Delta_{n,\lambda}=(1-\lambda)^{2}+\lambda^{2}\geq\frac{1}{2}$
for all $\lambda\in [0,1]$). For $n=1,2$, put
    $$
    f_{2}=x_{1}^{*}+2x_{2}^{*}+\frac{7}{3} x_{3}^{*}, \quad \textrm{and} \quad
    f_{1}=x_{1}^{*};
    $$
then $f_{2}(v_{3})=\frac{1}{3}$, $f_{2}(v_{2})=1$,
$f_{2}(v_{1})=1$, $f_{3}(v_{2})=-2$, $f_{1}(v_{2})=-1$,
$f_{1}(v_{1})=1$, and everything is fine (indeed,
$\Delta_{1,\lambda}=1$ and $\Delta_{2,\lambda}=(1-\lambda)^{2}+
\lambda^{2}\geq\frac{1}{2}$ for all $\lambda\in [0, 1]$).

Therefore, with these definitions, the linear system (3) has a
unique solution for $A_{n}, B_{n}$, which can be easily calculated
and estimated by Cramer's rule, of the form
    $$
    \begin{array}{ll}
    &A_{n}(y)=\frac{1}{\Delta_{n,\lambda}}\big(f_{n}(y)-
    \lambda [f_{n}(v_{n+1})-1]f_{n+1}(y)\big)\\
    &B_{n}(y)=\frac{1}{\Delta_{n,\lambda}}\big(f_{n+1}(y)-
    (1-\lambda)[f_{n+1}(v_{n})-1]f_{n}(y)\big).
    \end{array}
    $$
The linear forms $y\mapsto A_{n}(y)$, $y\mapsto B_{n}(y)$ satisfy
that $\|A_{n}\|\leq 144\geq \|B_{n}\|$ for all $n$, as is easily
checked. Now, by substituting $f_{n}(x)=A_{n}(y)$ and
$f_{n+1}(x)=B_{n}(y)$ in (2) we get the expression for the inverse
of $L_{n,\lambda}$, that is,
    $$
    x=L^{-1}_{n,\lambda}(y)=y-[(1-\lambda)A_{n}(y)(v_{n}-v_{1})+
    \lambda B_{n}(y)(v_{n+1}-v_{1})].  \eqno (4)
    $$
By taking into account that $\|A_{n}\|\leq 144\geq \|B_{n}\|$,
$\|f_{n}\|\leq 18$ and $\|v_{n}-v_{1}\|\leq 4$ for all $n$, one
can estimate that $1\leq\|L_{n,\lambda}\|\leq 73$ and
$1\leq\|L^{-1}_{n,\lambda}\|\leq 577$ for all $n\in\mathbb{N},
\lambda\in [0,1]$.

So let us define $\beta:[0,\infty)\longrightarrow\textrm{Isom}(X)$
by
    $$
    \beta(t)=\sum_{n=1}^{\infty}\theta_{n}(t)S_{n}. \eqno (5)
    $$
This path is well defined because the sum is locally finite; in
fact, from the definition of $\theta_n$ it is clear  that, for a
given $t_{0}\in[0,\infty)$ there exist some $\delta>0$ and
$N=N(t_{0})\in\mathbb{N}$ such that
$\beta(t)=\theta_{N}(t)S_{N}+\theta_{N+1}(t)S_{N+1}$ for all $t\in
(t_{0}-\delta, t_{0}+\delta)$, that is, $\beta$ is locally of the
form $\beta(t)=L_{n,\lambda(t)}$, where
$\lambda(t)=\theta_{n}(t)$. This implies that the $\beta(t)$ are
really linear isomorphisms and that the path is $C^{\infty}$
smooth.

On the other hand, the path
$\beta^{-1}(t)=[\beta(t)]^{-1}\in\textrm{Isom}(X)$ is
$C^{\infty}$ smooth as well, because it is the composition of our
path $\beta$ with the mapping
$\varphi:\textrm{Isom}(X)\longrightarrow\textrm{Isom}(X)$,
$\varphi(U)=U^{-1}$, which is $C^{\infty}$ smooth and whose
derivative is given by $\varphi'(U)(S)=-U^{-1}\circ S\circ
U^{-1}$ for every $S\in\mathcal{L}(X,X)$ (see \cite{Cartan},
theorem 5.4.3). This proves condition (1) of the lemma.

Next, by bearing in mind the local expression of $\beta$ and the
above estimations for $\|L_{n,\lambda}\|$ and
$\|L^{-1}_{n,\lambda}\|$, we deduce that
    $$
    1\leq\|\beta(t)\|\leq R\geq \|\beta^{-1}(t)\|\geq 1
    $$
for all $t\in [0,\infty)$, where $R\geq 577$ will be fixed later.
This shows condition (2). Now, if $t_{0}\in[0,\infty)$ and we
write $\beta(t)=\theta_{N}(t)S_{N}+\theta_{N+1}(t)S_{N+1}$ for
$t\in (t_{0}-\delta, t_{0}+\delta)$ as above, then it is clear
that $\beta'$ is locally of the form
    $$
    \beta'(t)=\theta'_{N}(t)S_{N}+\theta'_{N+1}(t)S_{N+1}
    $$
and therefore
    $$
    \|\beta'(t)\|\leq|\theta'_{N}(t)|\|S_{N}\|+|\theta'_{N+1}(t)|\|S_{N+1}\|\leq
    4(73+73)=584,
    $$
from which we get $\sup_{t\geq 0}\|\beta'(t)\|\leq 584\leq R$.
Moreover, we have
    $$
    (\beta^{-1})'(t)=-(\beta(t))^{-1}\circ \beta'(t)\circ (\beta(t))^{-1}
    $$
and therefore
    $$
    \|(\beta^{-1})'(t)\|\leq\|\beta(t)^{-1}\|^{2}\|\beta'(t)\|\leq
    (577)^{2}584,
    $$
from which $\sup_{t\geq 0}\|(\beta^{-1})'(t)\|\leq R$ and
condition (3) is satisfied as well provided we fix
$R=(577)^{2}584$.

Now let us define the path $p:[0,\infty)\longrightarrow X$ by
    $$
    p(t)=\int_{-\infty}^{t}\beta(s)(v_{1})ds=
    \int_{-\infty}^{t}\big(\sum_{n=1}^{\infty}\theta_{n}(s)S_{n}(v_{1})\big)ds.
    $$
It is clear that $p$ is a $C^{\infty}$ smooth path in $X$, and
$p'(t)=\beta(t)(v_{1})$ for all $t\geq 0$ (from which it follows
that $p$ is Lipschitz). Let us see that $p$ is bounded. For a
given $t>0$ there exists $N=N(t)\in\mathbb{N}$ so that
$N-1-\frac{1}{2}\leq t\leq N-\frac{1}{2}$ and therefore, taking
into account the definition of $\theta_n$ and the fact that
$S_{n}(v_{1})=v_{n}=x_{n}-x_{n-1}$ for all $n$, we have that
    \begin{eqnarray*}
    & & \|p(t)\|=\|\int^{t}_{-\infty}\sum_{n=1}^{\infty}\theta_{n}(s)S_{n}(v_{1})ds\|
    =\| \sum_{n=1}^{\infty}\big(\int^{t}_{-\infty}\theta_{n}(s)ds\big)v_{n}\|\\
    & &=\|\big(\int^{\infty}_{-\infty}\theta(s)ds\big)\sum_{n=1}^{N-1}v_{n}+
    \big(\int^{t}_{-\infty}\theta_{N}(s)ds\big)v_{N}\|\\
    & &\leq\big(\int^{\infty}_{-\infty}\theta(s)ds\big)\|\sum_{n=1}^{N-1}v_{n}\|
    +\big(\int^{\infty}_{-\infty}\theta(s)ds\big)\|v_{N}\|\\
    & &=\big(\int^{\infty}_{-\infty}\theta(s)ds\big)
    \big(\|x_{N-1}-x_{0}\|+\|x_{N}-x_{N-1}\|\big)
    \leq\frac{3}{2}\big(2+2)=6.
    \end{eqnarray*}
This shows that the image of $p$ is contained in the ball $B(0,
6)$ and $p$ is bounded. Let us also remark that
$2\geq\|v_{1}\|\geq\frac{x_{1}^{*}(x_{1}-x_{0})}{\|x_{1}^{*}\|}\geq\frac{1}{4}$.

Finally, let us check that $p$ satisfies the separation condition
(5). Let $0\leq t<r$ and take $N\in\mathbb{N}$ so that
$N-1-\frac{1}{2}<r\leq N-\frac{1}{2}$; then we have
    \begin{eqnarray*}
    & & p(r)-p(t)=\sum_{n=1}^{\infty}
    \big(\int_{t}^{r}\theta_{n}(s)ds\big)v_{n}=
    \sum_{n=1}^{\infty}\big(\int_{t}^{r}\theta_{n}(s)ds\big)(x_{n}-x_{n-1})=\\
    & &-\big(\int_{t}^{r}\theta_{1}(s)ds\big)x_{0}+
    \sum_{k=1}^{N-1}\big(\int_{t}^{r}\theta_{k}(s)ds-
    \int_{t}^{r}\theta_{k+1}(s)ds\big)x_{k} +
    \big(\int_{t}^{r}\theta_{N}(s)ds\big)x_{N}.
    \end{eqnarray*}
By observing that $\max\{1-s, 2s-1\}\geq\frac{1}{3}$ for all
$s\in\mathbb{R}$ and taking into account the definition of the
$\theta_n$, it is not difficult to see that
    $$
    \max\{\int_{t}^{r}\theta_{N}(s)ds,
    \int_{t}^{r}\theta_{N-1}(s)ds-\int_{t}^{r}\theta_{N}(s)ds\}\geq
    \min\{\frac{1}{3}|t-r|, a\}, \quad \eqno (6)
    $$
where $a=\int_{-\frac{1}{2}}^{0}\theta(s)ds>0$. Then, by
applying either $x_{N}^{*}$ or $x_{N-1}^{*}$ to the expression for
$p(r)-p(t)$ above, depending on which the maximum in (6) is, and
bearing in mind that $x_{n}^{*}(x_{k})=\delta_{n,k}$ and
$\|x_{n}^{*}\|\leq 4$ for all $n, k$, we get that
    $$
    \max\{x_{N}^{*}(p(r)-p(t)), x_{N-1}^{*}(p(r)-p(t))\}
    \geq\min\{\frac{1}{3}|t-r|, a\},
    $$
and it follows that $\|p(r)-p(t)\|\geq\min\{\frac{1}{12}|t-r|,
\frac{a}{4}\}$. This shows that if $R>0$ is large enough then
    $$
    \|p(r)-p(t)\|\geq\frac{1}{R}\min\{1, |t-r|\}
    $$
for all $t,r\geq 0$.

In order to get paths $\beta$ and $p$ and a vector $v$ with
properties (1)--(5) and such that $p$ is contained in the unit
ball, it is enough to multiply them all by $\frac{1}{6}$.
\end{proof}

\vspace{0.3cm}

\begin{center}
{\bf Proof of theorem \ref{Existence of twisted tubes}}.
\end{center}

Consider $X=H\oplus [z]=H\times\mathbb{R}$ and
$C_{\varepsilon}=\{x+tz\in X : \|x\|<\varepsilon, t>0\}$, where
$H=\textrm{Ker} z^{*}$ for some $z^{*}\in X^*$ with
$z^{*}(z)=\|z^*\|=\|z\|=1$, and $\varepsilon>0$ is to be fixed
later. Let $\beta$ and $p$ be the paths from lemma \ref{Path beta
of linear isomorphisms}. There is no loss of generality if we
assume that $v\in [z]$, $z^{*}(v)\geq\frac{1}{K}$. Let us define
$\pi:C_{\varepsilon}\longrightarrow X$ by
    $$
    \pi(x,t)=\beta(t)(x)+p(t).
    $$
It is clear that $\pi$ is $C^{\infty}$ smooth and has a bounded
derivative. We are going to show that $\pi$ is a diffeomorphim
onto its image, $T_{\varepsilon}$, and
$\pi^{-1}:T_{\varepsilon}\longrightarrow C_{\varepsilon}$ has a
bounded derivative as well. To this end let us define the path
$\alpha:[0,\infty)\longrightarrow X^*$ by
    $$
    \alpha(t)=f_{t}=z^{*}\circ \beta^{-1}(t).
    $$
This is a $C^{\infty}$ smooth and Lipschitz path in $X^*$, and
$\alpha(t)=f_{t}$ satisfies that $\textrm{Ker}f_{t}=\beta(t)(H)$.
It is clear from this definition and the properties of $\beta$ and
$p$ that
\begin{itemize}
\item [{(i)}] $\|\alpha'(t)\|\leq K$, and
\item [{(ii)}] $\alpha(t)(p'(t))=z^{*}(v)\geq\frac{1}{K}$
\end{itemize}
for all $t\geq 0$. Now, for a fixed (but arbitrary) $y\in
T_{\varepsilon}=\pi(C_{\varepsilon})$, let us introduce the
auxiliary function $F=F_{y}:[0, \infty)\longrightarrow\mathbb{R}$
defined by
    $$
    F(t)=\alpha(t)(y-p(t)).
    $$
We have that
    \begin{eqnarray*}
    & & F_{y}'(r)=\alpha'(r)(y-p(r))-\alpha(r)(p'(r))\\
    & &\leq\|\alpha'(r)\| \|y-p(r)\| -\alpha(r)(p'(r))\\
    & &\leq K\|y-p(r)\| -\frac{1}{K}
    \end{eqnarray*}
for all $r\geq 0$. If we choose $\varepsilon>0$ smaller than
$\frac{1}{6K^{5}}$ this implies that $\pi$ is a $C^{\infty}$
diffeomorphism onto its image.

Indeed, let us first see that $\pi$ is an injection. Assume that
$y=\pi(x,t)=\pi(w,s)$ for some $(x,t), (w,s)\in C_{\varepsilon}$.
Then we have $y-p(t)=\beta(t)(x)$ and $y-p(s)=\beta(s)(w)$, so
that $x=\beta^{-1}(t)(y-p(t))$ and $w=\beta^{-1}(s)(y-p(s))$,
and, in order to conclude that $(x,t)=(w,s)$, it is enough to see
that $t=s$. Note that $\beta(t)(x)-\beta(s)(w)=p(s)-p(t)$ and
therefore, by (5) of lemma \ref{Path beta of linear
isomorphisms},
    \begin{eqnarray*}
    & & \frac{1}{K}\min\{1, |t-s|\}\leq\|p(s)-p(t)\|=
    \|\beta(t)(x)-\beta(s)(w)\|\\
    & &\leq\|\beta(t)(x)\|+ \|\beta(s)(w)\|\leq K(\|x\|+\|w\|)\leq
    2K\varepsilon\leq\frac{1}{3K^{4}},
    \end{eqnarray*}
so that $|t-s|\leq 2K^{2}\varepsilon\leq\frac{1}{3K^{3}}$. Now,
since $p$ and $\beta$ are both $K$-Lipschitz, for every $r\in
[t,s]$ we have that
    \begin{eqnarray*}
    & &\|y-p(r)\|\leq\|y-p(t)\|+\|p(t)-p(r)\|=\|\beta(t)(x)\|+\|p(t)-p(r)\|\\
    & &\leq K\|x\|+K|t-r|\leq K\varepsilon +2K^{3}\varepsilon\leq
    3K^{3}\varepsilon.
    \end{eqnarray*}
By combining this with the above estimation for $F_{y}'(r)$ we get

    $$
    F_{y}'(r)\leq K\|y-p(r)\| -\frac{1}{K}\leq 3K^{4}\varepsilon
    -\frac{1}{K}\leq -\frac{1}{2K} \quad \eqno (7)
    $$
for every $r\in [t,s]$. Now suppose that $t\neq s$. Then, since
$x=\beta^{-1}(t)(y-p(t))$ and $w=\beta^{-1}(s)(y-p(s))$ are both
in $H$ we have that $0=z^{*}(x)=z^{*}(w)=F_{y}(t)=F_{y}(s)$, so
that, by the classic Rolle's theorem, there should exist some
$r\in (t,s)$ with $F_{y}'(r)=0$. But this contradicts (7).
Therefore $t=s$ and $\pi$ is an injection.

If, for a given $y\in\pi(C_{\varepsilon})$, we denote by $t(y)$
the unique $t=t(y)$ such that $y=\pi(\beta^{-1}(t)(y-p(t)), t)$
then it is clear that the inverse
$\pi^{-1}:T_{\varepsilon}\longrightarrow C_{\varepsilon}$ is
defined by
    $$
    \pi^{-1}(y)=(\beta^{-1}(t(y))(y-p(t(y))), t(y)). \quad \eqno (8)
    $$
For each $y$ the number $t(y)$ is uniquely determined by the
equation
    $$
    G(y, t):=F_{y}(t)=0,
    $$
and the argument above shows that
    $$
    \frac{\partial G}{\partial t}(y,t)=F_{y}'(t)\leq
    -\frac{1}{2K} \quad \eqno (9)
    $$
for every $y\in T_{\varepsilon}$ and $t$ in a neighbourhood of
$t(y)$. Then, according to the implicit function theorem we get
that the function $y\mapsto t(y)$ is $C^{\infty}$ smooth.
Furthermore, we have that
    $$
    t'(y)=\frac{-\frac{\partial G}{\partial y}(y,t(y))}
    {\frac{\partial G}{\partial t}(y,t(y))}=
    \frac{-z^{*}\circ\beta^{-1}(t(y))}{F_{y}'(t(y))},
    $$
and therefore, according to the above estimations,
    $$
    \|t'(y)\|\leq
    \|z^{*}\circ\beta^{-1}(t(y)) \|\frac{1}{| F_{y}'(t(y))|}\leq
    2K^{2},
    $$
which shows that $y\mapsto t(y)$ has a bounded derivative as well.
Then it is clear that $\pi^{-1}$ is $C^{\infty}$ and has a bounded
derivative (all the functions involved in (8) have been proved to
have bounded derivatives). This concludes the proof of theorem
\ref{Existence of twisted tubes}.

\vspace{0.2cm}

We will finish this section with a simple alternative proof of
the failure of Rolle's theorem in the
non-Lipschitz case.

\begin{rem}
{\em If we drop the Lipschitz condition from
the statement of theorem \ref{Rolle's theorem fails}, a much
simpler proof based on the same idea is available. Let us make a
sketch of this proof.}
\end{rem}
Consider the decomposition $X=H\times\mathbb{R}$ and pick a
non-negative $C^{p}$ smooth bump function $\varphi$ on $H$ whose
support is contained on the ball $B_{H}(0,1/16)$. First, we
construct a $C^{\infty}$ smooth path $q:[0,\infty)\longrightarrow
B_H$, where $B_H$ stands for the unit ball of the hyperplane $H$,
with the property that $q$ has no accumulation points at the
infinity, that is, $\lim_{n\to\infty}q(t_{n})$ does not exist for
any $(t_n)$ going to $\infty$. This can easily be done by having
$q$ lost in the infinitely many dimensions of $H$. For instance,
take a biorthogonal sequence $\{x_{n}, x_{n}^{*}\}\subseteq
H\times H^*$ so that $\|x_{n}\|=1$ and $\|x_{n}^{*}\|\leq 4$, and
consider a $C^{\infty}$ function $\theta:\mathbb{R}\longrightarrow
[0,1]$ so that $\textrm{supp}\theta\subseteq [-1,1]$,
$\theta(0)=1$, $\theta'(t)<0$ for $t\in (0,1)$, and
$\theta(t-1)=1-\theta(t)$ for $t\in [0,1]$. The path $q$ may be
defined as
    $$
    q(t)=\sum_{n=1}^{\infty}\theta(t-n+1)x_{n}
    $$
for $t\geq 0$. Now we reparametrize $q$ and define
$p:[0,1)\longrightarrow B_H$ by
    $$
    p(t)=q\big(\frac{t}{1-t}\big).
    $$
Let $\alpha:\mathbb{R}\longrightarrow [0,1]$ be a $C^{\infty}$
smooth function so that $\alpha(t)=0$ for all $t\leq 0$, and
$\alpha'(t)>0$ for all $t>0$. Then the function
$g:X=H\times\mathbb{R}\longrightarrow\mathbb{R}$ defined by
    $$
    g(x,t)=\left\{
    \begin{array}{ll}
    &\varphi(x-p(t))\alpha(t) \quad \textrm{if} \quad t\in [0,1);\\
    &0 \quad \textrm{otherwise}
    \end{array}\right.
    $$
is a $C^p$ smooth bump function which does not satisfy Rolle's
theorem. Indeed, it is easy to see that
    $$
    g'(x,t)(p'(t),1)=\varphi(x-p(t))\alpha'(t)>0,
    $$
and in particular $g'(x,t)\neq 0$, for all $(x,t)$ in the interior
of the support of $g$.


\section[Killing singularities]{Killing singularities. The failure of Brouwer's
fixed point theorem in infinite dimensions.}

Do not be afraid, this section does not contain any totalitarian propaganda.
Here we will present two applications of theorem
\ref{Characterization by diffeomorphisms deleting twisted tubes},
both of which have in common the following principle: if you have
a mapping with a single singular point or an isolated set of
singularities that bother you, you can just kill them by
composing your map with some deleting diffeomorphisms. In this way
you obtain a new map which is as close as you want to the old
one but does not have the adverse properties created by the
singular points you eliminate.

For instance, if you want
a smooth bump function $g$ which does not satisfy Rolle's theorem and
whose support is a smooth starlike body $A$, by composing the
Minkowski functional of this body with a real bump function you
get a function $h$ whose support is $A$ and whose derivative
vanishes only at the origin and outside $A$; then, by composing $h$
with a diffeomorphism $f$ which extracts a small set containing the origin
and which restricts to the identity outside $A$, you get a map $g$
with the required properties.

On the other hand, suppose you want a smooth retraction $r$ from a
bounded starlike body $A$ of a Banach space $X$ onto its boundary
$\partial A$. This is impossible if $X$ is finite-dimensional, but
otherwise you can use the following trick: it is trivial
that there is a smooth retraction $h$ from $A\setminus\{0\}$ onto $\partial A$;
then take a diffeomorphism $f$
which removes from $X$ a small subset containing the origin and
restricts to the identity outside $A$. The composition $r=h\circ f$
gives the required retraction.

Let us formalize these ideas and comment on the results that they
provide.

\vspace{0.2cm}

\begin{center}
{\small{\bf The support of the bumps that violate Rolle's theorem}}.
\end{center}

The bump function constructed in the proof of theorem \ref{Rolle's
theorem fails} has a weird support, namely a twisted tube. Some
readers might judge this fact rather unpleasant and wonder
whether it is possible to construct a bump function which does
not satisfy Rolle's theorem and whose support is a nicer set,
such as a ball or a starlike body. To comfort those readers let
us first recall that in infinite dimensions there is no topological
difference between a tube (whether it is twisted or not) and a
ball or a starlike body, as theorem 2.4 in \cite{Adev2} shows.
Furthermore, as we said above, theorem
\ref{Characterization by diffeomorphisms deleting twisted tubes}
allows us to show that
for a given $C^p$ smooth bounded starlike body $A$ in an infinite-dimensional
Banach space $X$, it is always possible to construct a $C^p$ smooth
bump function on $X$ which does
not satisfy Rolle's theorem and whose support is precisely the
body $A$.

\begin{thm}
Let $X$ be an infinite-dimensional Banach space with a $C^p$
smooth bounded starlike body $A$. Then there exists a $C^p$ smooth
bump function $g$ on $X$ whose support is precisely the body $A$,
and with the property that $g'(x)\neq 0$ for all $x$ in the
interior of $A$ (that is, $g$ does not satisfy Rolle's theorem).
\end{thm}
\begin{proof}
Let $q_A$ be the Minkowski functional of $A$.
We may assume that $B_{X}\subseteq A$. By theorem
\ref{Characterization by diffeomorphisms deleting twisted tubes}
there is a closed subset $D$ of $A$ and a $C^p$ diffeomorphism
$f:X\longrightarrow X\setminus D$ which is the identity outside
$A$. It can be assumed that the origin belongs to
$D$. Then the function $h:X\longrightarrow\mathbb{R}$ defined by
    $$
    h(x)=q_{A}(f(x))
    $$
is $C^p$ smooth on $X$, restricts to the gauge $q_A$ outside $A$,
and has the remarkable property that $h'(x)\neq 0$ for all $x\in X$
(indeed, $h'(x)=q_{A}'(f(x))\circ f'(x)$ is non-zero everywhere
because $q_{A}'(y)\neq 0$ whenever $y\neq 0$, $0\notin f(X)$, and
$f'(x)$ is a linear isomorphism at each point $x$).

Now, take a $C^\infty$
real function $\theta:\mathbb{R}\longrightarrow [0, 1]$ such that
$\theta(t)>0$ for $t\in (-1, 1)$, $\theta=0$ outside $[-1, 1]$,
$\theta(t)=\theta(-t)$, $\theta(0)=1$, and $\theta'(t)<0$ for all
$t\in (0, 1)$. Then, if we define $g:X\longrightarrow\mathbb{R}$ by
    $$
    g(x)=\theta(h(x)),
    $$
it is immediately checked that $g$ is a $C^p$ smooth bump on $X$ which
does not satisfy Rolle's theorem and whose support is precisely
the body $A$.
\end{proof}

\begin{center}
{\small{\bf The failure of Brouwer's fixed point theorem in infinite dimensions}}.
\end{center}

The celebrated Brouwer's fixed point theorem tells us that every
continuous self-map of the unit ball of a finite-dimensional
normed space admits a fixed point. This is the same as saying
that there is no continuous retraction from the unit ball onto
the unit sphere, or that the unit sphere is not contractible (the
identity map on the sphere is not homotopic to a constant map).
The {\em Scottish book} (see \cite{Mauldin}) contains the
following question, asked in 1935 by S. Ulam. Can one transform
continuously the solid sphere of a Hilbert space into its
boundary such the transformation should be the identity on the
boundary of the ball? In other words, is the unit sphere of a
Hilbert space a retract of its unit ball? The first answer to
this question is commonly attributed to S. Kakutani
\cite{Kakutani}, who solved the problem by exhibiting several
examples of continuous self-mappings  of the unit ball of the
Hilbert space without fixed points. Thus, none of the above forms
of Brouwer's fixed point theorem remains valid in infinite
dimensions.

A very nice solution to the retraction problem, and one which has
the advantage of holding in arbitrary infinite-dimensional Banach
spaces, was given by the pioneering results of Klee's on
topological negligibility of points and compacta \cite{Klee2,
Klee1}: for every infinite-dimensional Banach space $X$ there
always exists a homeomorphism $h:X\longrightarrow
X\setminus\{0\}$ so that $h$ restricts to the identity outside
the unit ball $B_X$. The required retraction of $B_X$ onto the
unit sphere $S_X$ is then given by $R(x)=h(x)/\|h(x)\|$ for $x\in
B_X$. By taking into account the subsequent progress on
topological negligibility of subsets made by C. Bessaga, T.
Dobrowolski and the first-named author among others (see
\cite{Be, BP, Do, Do2, A, Ado}), this mapping $h$ may even be
assumed $C^p$ smooth provided that the sphere $S_X$ is $C^p$
smooth. Thus, there are very regular retractions that provide an
answer to Ulam's question.

In \cite{Nowak} B. Nowak showed that for several
infinite-dimensional Banach spaces Brouwer's theorem fails even
for Lipschitz mappings (that is, under the strongest
uniform-continuity condition), and in \cite{BS} Y. Benyamini and
Y. Sternfeld generalized Nowak's result for all
infinite-dimensional normed spaces, establishing that for every
infinite-dimensional space $(X, \|\cdot\|)$ there exists a
Lipschitz retraction from the unit ball $B_{X}$ onto the sphere
$S_{X}$, and that $S_X$ is Lipschitz contractible.

More recently, M. Cepedello and the first-named author showed that
these results hold for the smooth category as well (see
\cite{AC}). In fact they proved that for every
infinite-dimensional Banach space with a $C^p$ (Lipschitz) bounded
starlike body $A$ (where $p=0, 1, 2, \dots, \infty$), there is a
$C^p$ Lipschitz retraction of $A$ onto its boundary $\partial A$,
the boundary $\partial A$ is $C^p$ (Lipschitz) contractible, and there
is a $C^p$ smooth (Lipschitz) mapping $f:A\longrightarrow A$ such that
$f$ has no (approximate) fixed points.

The proof of these results in the general case is somewhat involved,
but if we drop the Lipschitz condition then the fact that Brouwer's
theorem is false in infinite dimensions even for smooth self-mappings
of balls or starlike bodies  is a trivial consequence of theorem
\ref{Characterization by diffeomorphisms deleting twisted tubes}.

\begin{cor}[Azagra--Cepedello]\label{Retractions of starlike bodies}
Let $X$ be an infinite-dimensional Banach space and let $A$ be a
$C^p$ smooth bounded starlike body. Then:
\begin{itemize}
\item [{(1)}] The boundary $\partial A$ is $C^p$ contractible.
\item [{(2)}] There is a $C^p$ smooth retraction from $A$ onto
$\partial A$.
\item [{(3)}] There exists a $C^p$ smooth mapping $\varphi:A\longrightarrow A$
without fixed points.
\end{itemize}
\end{cor}
\begin{proof}
Let $f:X\longrightarrow X\setminus D$ be the
diffeomorphism from theorem
\ref{Characterization by diffeomorphisms deleting twisted tubes}.
We may assume that the origin belongs to the deleted set $D$ and
that $B_{X}\subseteq A$, so that
$f$ restricts to the identity outside $A$. Then the formula
    $$
    R(x)=\frac{f(x)}{q_{A}(f(x))},
    $$
where $q_A$ is the Minkowski functional of $A$, defines a $C^p$
smooth retraction from $A$ onto the boundary $\partial A$. This
proves (2).

Once we have such a retraction it is easy to prove parts (1) and (3):
the formula $\varphi(x)=-R(x)$ defines a $C^p$ smooth self-mapping
of $A$ without fixed points. On the other hand, if we pick a non-decreasing
$C^{\infty}$ function $\zeta:\mathbb{R}\longrightarrow\mathbb{R}$
so that $\zeta(t)=0$ for $t\leq\frac{1}{4}$ and $\zeta(t)=1$ for
$t\geq\frac{3}{4}$, then the formula
    $$
    H(t,x)=R((1-\zeta(t))x),
    $$
for $t\in[0, 1]$, $x\in\partial A$, defines a $C^p$
homotopy joining the identity to a constant on $\partial A$, that
is, $H$ contracts the pseudosphere $\partial A$ to a point.
\end{proof}

\section{How small can the set of gradients of a bump be?}

If $b$ is a smooth bump function on a Banach space $X$ it is
natural to ask how large or how small the cone generated by the
set of gradients $b'(X)$ can be. In general, as a consequence of
Ekeland's variational principle, one has that the cone
$\mathcal{C}(b)=\{\lambda b'(x) : \lambda\geq 0, x\in X\}$
is norm-dense in the dual space $X^*$ (see \cite{DGZ}, pag. 58,
proposition 5.2).

In \cite{Adev2} a study was initiated on the
topological properties of the set of derivatives of smooth
functions. Among other results it was proved that an
infinite-dimensional separable Banach space has a
$C^1$ smooth bump function (resp. is Asplund) if and only if there exists
another $C^1$ smooth bump function $b$ on $X$ with the property that
$b'(X)=X^{*}$. This answers the question as to how large can
the cone $\mathcal{C}(b)$ be. But what is the smallest possible
size of $\mathcal{C}(b)$?

To begin with, by using theorem \ref{Rolle's theorem fails}, one
can easily construct smooth bump functions
whose sets of gradients lack not only the point zero, but any
pre-set finite-dimensional linear subspace of the dual.

\begin{cor}
Let $X$ be an infinite-dimensional Banach space and $W$ a
finite-dimensional subspace of $X^*$. The following statements are
equivalent.
\begin{itemize}
\item [{(1)}] $X$ has a $C^p$ smooth (Lipschitz) bump function.
\item [{(2)}] $X$ has a $C^p$ smooth (Lipschitz) bump function $f$
satisfying that $\{\lambda f'(x) : x\in X,
\lambda\in\mathbb{R}\}\cap W=\{0\}$. Moreover, $\{f'(x): f(x)\neq
0\}\cap W=\emptyset$.
\end{itemize}
\end{cor}
\begin{proof}
We can write $X=Y\oplus Z$, where $Y=\cap_{w^{*}\in
W}\textrm{Ker}w^{*}$ and $\textrm{dim}Z=\textrm{dim}W$ is finite.
Pick a $C^p$ smooth (Lipschitz) bump function $\varphi:
Y\longrightarrow\mathbb{R}$ such that $\varphi$ does not satisfy
Rolle's theorem, and let $\theta$ be a $C^{\infty}$ smooth
Lipschitz bump function on $Z$ so that $\theta'(z)=0$ whenever
$\theta(z)=0$. Then the function $f:X=Y\oplus
Z\longrightarrow\mathbb{R}$ defined by
$f(y,z)=\varphi(y)\theta(z)$ is a $C^p$ smooth (Lipschitz) bump
which satisfies $\{f'(x): f(x)\neq 0\}\cap W=\emptyset$. Indeed,
we have
    $$
    f'(y,z)=\big(\theta(z)\varphi'(y), \varphi(y)\theta'(z)\big)\in
    X^{*}=Y^{*}\oplus Z^{*}=Y^{*}\oplus W
    $$
and, since $\theta(z)\varphi'(y)\neq 0$ whenever
$\varphi(y)\theta'(z)\neq 0$, it is clear that $f'(y,z)\notin W$
unless $f'(y,z)=0$.
\end{proof}

We also have the following

\begin{cor}
Let $X$ be an infinite-dimensional Banach space such that
$X=X_{1}\oplus X_{2}$, where $X_{1}$ and $X_{2}$ are both
infinite-dimensional. The following statements are equivalent.
\begin{itemize}
\item [{(1)}] $X$ has a $C^p$ smooth (Lipschitz) bump function.
\item [{(2)}] $X$ has a $C^p$ smooth (Lipschitz) bump function $f$
satisfying that $\{\lambda f'(x) : x\in X,
\lambda\in\mathbb{R}\}\cap (X_{1}^{*}\cup X_{2}^{*})=\{0\}$.
Moreover, $\{f'(x): f(x)\neq 0\}\cap (X_{1}^{*}\cup
X_{2}^{*})=\emptyset$.
\end{itemize}
\end{cor}
\begin{proof}
The proof is similar to that of the preceding corollary. Pick
$\varphi_{1}$ and $\varphi_{2}$ smooth (Lipschitz) bump functions on
$X_{1}$ and $X_{2}$, respectively, so that $\varphi_{1}$ and
$\varphi_{2}$ do not satisfy Rolle's theorem. Then the function
$f:X=X_{1}\oplus X_{2}\longrightarrow\mathbb{R}$ defined by
$f(x_{1},x_{2})=\varphi_{1}(x_{1})\varphi_{2}(x_{2})$ is a smooth
(Lipschitz) bump which satisfies $\{f'(x): f(x)\neq 0\}\cap
(X_{1}^{*}\cup X_{2}^{*})=\emptyset$.
\end{proof}

If we restrict the scope of our search to classic Banach spaces,
much stronger results are available.
On the one hand, if $X=c_{0}$ the size of $\mathcal{C}(b)$ can be
really small. Indeed, as a consequence of P. H\'ajek's work \cite{Hajek}
on smooth functions on $c_0$ we know that if $b$ is $C^1$ smooth
with a locally uniformly continuous derivative (note that there are
bump functions with this property in $c_0$), then $b'(X)$ is
contained in a countable union of compact sets in $X^*$ (and in
particular has empty interior).  On the other hand, if $X$ is
non-reflexive and has a separable dual, there are bumps $b$ on $X$ so that $\mathcal{C}(b)$
has empty interior, as it was shown in \cite{Adev2}.

In the reflexive case, however, the problem is far from being
settled. To begin with, the cone $\mathcal{C}(b)$ cannot be very
small, since it is going to be a residual subset of the dual $X^*$.
Indeed, as a straightforward consequence of Stegall's variational
principle, for every Banach space $X$ having the Radon-Nikodym Property
(RNP) it is easy to see that $\mathcal{C}(b)$ is a
residual set in $X^*$. Therefore, for
infinite-dimensional Banach spaces $X$ enjoying RNP one can hardly
expect a better answer to the above question than the following
one: there are smooth bumps $b$ on $X$ such that the cone
$\mathcal{C}(b)$ has empty interior in $X^*$. In the case of the Hilbert space
$X=\ell_{2}$ we next show that this indeed happens.

\begin{thm}\label{eledos}
In the Hilbert space $\ell_2$ the following holds:
\begin{itemize}
\item [{(1)}] The usual norm $||\cdot||_{2}$ can be uniformly approximated
by  $C^1$  Lipschitz   functions $\psi$ (with Lipschitz
derivative) so that the cones $\mathcal{C}(\psi)$
generated by the sets of derivatives of $\psi$
have empty interior, and $\psi'(x)\neq 0$ for all
$x\in\ell_{2}$.
\item [{(2)}] There is a $C^1$ Lipschitz bump function $b$
(with Lipschitz derivative) on $\ell_2$ satisfying that the cone
$\mathcal{C}(b)$
generated by its set of derivatives $b'(\ell_{2})$ has empty
interior, and $b'(x)\neq 0$ for every $x$ in the interior of its support.
\end{itemize}
\end{thm}
\begin{proof}
To save notation, let us just write $||\cdot||$ when
referring  to the usual norm in $\ell_2$. We will make use of the following
restatement of a striking result due to S. A. Shkarin (see
\cite{Shk}).

\begin{thm}[Shkarin]
There is a $C^\infty$ diffeomorphism $\varphi$  from $\ell_2$ onto
$\ell_2\setminus\{0\}$ such that all the derivatives $\varphi^{(n)}$ are uniformly
continuous on $\ell_2$, and $\varphi(x)=x$  for $||x|| \ge 1$.
\end{thm}

Let us consider, for $0<\varepsilon<1$, the difeomorphism
$\varphi_\varepsilon:\ell_2 \longrightarrow \ell_2\setminus\{0\}$,
$\varphi_\varepsilon(x)=\varepsilon \varphi(x/\varepsilon)$, and
the function $U=U_{\varepsilon}:\ell_2 \longrightarrow \mathbb{R}$ defined by
$U(x)=\varepsilon^2 + ||\varphi_\varepsilon(x)||^2$. Then $U$
satisfies the following properties:
\begin{itemize}
\item[{(i)}] $U$ is $C^{\infty}$ smooth.
\item[{(ii)}] $||x||^2\le U(x)\le 2\varepsilon^2 +||x||^2$ and $\varepsilon^2 \le U(x)$, for
every $x\in \ell_2$.
\item[{(iii)}] $U(x)=\varepsilon^2 + ||x||^2$, for every $x\in
\ell_2$, $||x||\ge \varepsilon$.
\item[{(iv)}] $U'(x)\not=0$ for every  $x\in \ell_2$.
\item[{(v)}] $U$ is Lipschitz in bounded sets and $U'$ is Lipschitz.
\end{itemize}
Now, we define the functions $U_n:\ell_2\longrightarrow \mathbb{R}$ by
$U_n(x)=\frac1{2^{2n}}U(2^nx)$, whenever $x\in \ell^2$. We will
identify
$\ell_2$ with the infinite sum $\sum_2\ell_2 \equiv
\ell_2\oplus_2\ell_2\oplus_2 \ell_2 \cdots$, where an element
$x=(x_n)$ belongs to $\sum_2 \ell_2$ if and only if every $x_n$ is
in $\ell_2$ and $\sum_n ||x_n||^2<\infty$, being $||x||^2=\sum_n
||x_n||^2$. Then, we define the function $f:\sum_2 \ell_2 \longrightarrow
\mathbb{R}$ by
    $$
    f(x)=\sum_n U_n(x_n), \text{ where } x=(x_n)_n.
    $$
First, note that $f$ is well-defined, since condition (ii)
implies that, whenever $x=(x_n)\in \sum_2\ell_2$,

\begin{eqnarray*}
& & 0<f(x)=\sum_n \frac1{2^{2n}}U(2^nx_n)\le \sum_n \frac1{2^{2n}}
(2\varepsilon^2+||2^nx_n||^2)\\
& &=\sum_n (\frac{2\varepsilon^2}{2^{2n}} +||x_n||^2)<\infty.
\hspace{5cm} (4.1)
\end{eqnarray*}

On the one hand, note that, if $U'$ has Lipschitz constant less
than or equal to
$M$ then $U'_n$ is also Lipschitz with constant less than or equal to $M$,
since for $x$ and $y$ in $\ell_2$ we have
     $$
     ||U'_n(x)-U'_n(y)||=\frac1{2^n}||U'(2^nx)-U'(2^ny)||\le
     \frac1{2^n}M 2^n||x-y||.
     $$
This implies that, if $x=(x_n)\in \sum_2
\ell_2$, the functionals $U'_n(x_n)\in \ell_2$ satisfy that
$(U'_n(x_n))_n \in \sum_2 \ell_2$. Indeed, we have
$||U'_n(x_n)-U'_n(0)||\le M ||x_n||$, and therefore
$\sum_n||U'_n(x_n)-U'_n(0)||^2 <\infty$. Also,
$(U'_n(0))=(\frac1{2^n}U'(0))\in (\sum_2 \ell_2)^*\equiv \sum_2
\ell_2$, and then we get that $T(x)=(U_n'(x_n))$  also belongs to
$\sum_2 \ell_2$.

Let us now prove that $f$ is $C^1$ smooth.  For every $x=(x_n)$ and
$h=(h_n)$ in $\sum_2 \ell_2$, we can estimate
\begin{eqnarray*}
& & |f(x+h)-f(x)-T(x)(h)|\le\sum_n
|U_n(x_n+h_n)-U_n(x_n)-U'_n(x_n)(h_n)| \\ & & \le \sum_n
|U'_n(x_n+t_nh_n)(h_n)-U'_n(x_n)(h_n)| \qquad \text{ (for some }
0\le t_n\le 1) \\ & & \le M \sum_n ||h_n||^2=M||h||^2.
\end{eqnarray*}
Therefore $f$ is Fr\'echet differentiable and
$f'(x)=(\frac1{2^n}U'(2^n x_n))$. Moreover, $f'$ is Lipschitz
since $||f'(x)-f'(y)||^2=\sum_n ||U'_n(x_n)-U'_n(y_n)||^2\le M^2
\sum_n ||x_n-y_n||^2=M^2 ||x-y||^2$. This implies, in particular,
that $f$ is Lipschitz on bounded sets.

Let us check that $f=f_{\varepsilon}$ uniformly approximates $||\cdot||^2$.
Indeed, from condition (ii) on $U$ and (4.1), we have
that, for every $x=(x_n)\in\sum_2\ell_2$,
     $$
     \max\{\frac 13 \varepsilon^2 ,||x||^2\}\le f(x)\le \frac 23
     \varepsilon^2 +||x||^2,  \eqno(4.2)
     $$
and then,
     $$
     0\le f(x)-||x||^2\le \frac 23 \varepsilon^2.  \eqno(4.3)
     $$
In order to obtain functions which approximate the
norm uniformly in $\ell_2$ let us consider
$\psi=\psi_{\varepsilon}=\sqrt{f_{\varepsilon}}$.
According to inequalities
(4.2) and (4.3) we have that
     $$0\le
     \psi-||x||\le\frac{2\varepsilon^2}{3(\psi+||x||)} \le
     \frac{2}{\sqrt{3}}\varepsilon
     $$
for any $x\in \sum_2 \ell_2$.

Let us check that $\psi'$ is  bounded. By equation (4.2)
we have,  for any $x\in \sum_2\ell_2$,
     $$||\psi'(x)||=\frac {||f'(x)||}{2\psi(x)}
     \le \frac {||f'(x)-f'(0)||}{2\psi(x)} +\frac {||f'(0)||}{2\psi(x)}\le
     \frac{M}{2} +
     \frac{\sqrt{3}}{2\varepsilon}||f'(0)||.
     $$
Consequently, $\psi$ is Lipschitz with Lipschitz constant, say N.
In a similar way, we obtain that $\psi'$ is Lipschitz, since for any
$x$, $y$ in $\sum_2\ell_2$,
\begin{eqnarray*}
& & ||\psi'(x)-\psi'(y)||=||\frac {f'(x)-f'(y)}{2\psi(x)}+\frac{f'(y)}{2}
\big( \frac{1}{\psi(x)}-\frac{1}{\psi(y)}\big) || \\
& & \le \frac{1}{2} \frac{||f'(x)-f'(y)||}{\psi(x)}+
\frac{||\psi(y)-\psi(x)||}{\psi(x)}\frac{||f'(y)||}{2\psi(y)} \\
& & \le
\frac{\sqrt{3}M}{2\varepsilon}||x-y|| +\frac{\sqrt{3}N^2}{\varepsilon}\|x-y\|.
\end{eqnarray*}

\medskip

Finally, note that the set $\{\lambda f'(x)=\lambda (U'_n(x_n)):
\,\,x=(x_n)\in \sum_2\ell_2, \lambda>0\}$ is contained in $\{z=(z_n)\in \sum_2\ell_2
:\,\, z_n\not=0 \, \text{ for every }n \in \mathbb{N}\}$, which  has empty
interior in $\sum_2\ell_2$. This concludes the proof of (1).

\medskip

In order to   prove (2), we consider a $C^\infty$  function
$\theta:\mathbb{R}^{+} \longrightarrow \mathbb{R}$, $\theta'(t)<0$ for
$t\in(0,1)$, and supp$\, \theta=(0,1]$. Then, we can define a
required bump function as the composition $b(x)=\theta(f(x))$.
Indeed, on the one hand, $f(0)\le \frac23\varepsilon^2<1$ and
therefore
$b(0)>0$. On the other hand, $f(x)\ge||x||^2\ge 1$, whenever
$||x||\ge 1$, and hence $b(x)=0$ for $||x||\ge 1$. The bump function
$b$ is clearly Lipschitz with Lipschitz derivative since $\theta$,
  $\theta'$ and $f'$ are  Lipschitz and $f$ is Lipschitz on bounded sets.
\end{proof}

It is clear that the proof of the preceding theorem could only be adapted for
superreflexive Banach spaces $X$ which admit a decomposition
as an infinite sum of subspaces isomorphic to $X$, and therefore
the problem whether in a separable reflexive space there can be bumps
whose sets of gradients have empty interior remains open in the general case, though the
information that has already been gathered seems to point to a
final positive solution.

\begin{op}
{\em Let $X$ be a (separable reflexive) infinite-dimensional Banach space which admits
a $C^p$ smooth (Lipschitz) bump function. Is there another $C^p$
smooth (Lipschitz) bump function $f$ on $X$ so that the cone
$\mathcal{C}(f)=\{\lambda f'(x) :\lambda\geq 0, x\in X\}$ has empty interior in
$X^{*}?$}
\end{op}

The proof of theorem \ref{eledos} naturally raises the following
question.

\begin{op}
{\em Can every equivalent norm on $\ell_2$ be uniformly approximated
by $C^1$ Lipschitz functions satisfying that the cone generated by their derivatives has
empty interior? Furthermore, is it possible to approximate squared
equivalent norms in $\ell_{2}$ uniformly on bounded sets by real-analytic functions
(resp. polynomials) $\psi$ such that the cones generated by the sets $\psi'(\ell_{2})$
have empty interior?}
\end{op}


\vspace{0.7cm}
\begin{center}
{\bf Acknowledgements}
\end{center}
This research was carried out during a postdoctoral stay of the
authors in the Equipe d'Analyse de l'Universit\'e Pierre et Marie
Curie, Paris 6. The authors are indebted to the Equipe d'Analyse
and very especially to Gilles Godefroy for their kind hospitality
and generous advice.


\vspace{3mm}

\noindent Departamento de An\'alisis Matem\'atico, Facultad de
Ciencias Matem\'aticas, Universidad Complutense, 28040 MADRID,
SPAIN.

\noindent Equipe d'Analyse, Universit\'e Pierre et Marie
Curie--Paris 6. 4, place Jussieu, 75005 PARIS, FRANCE.

\noindent {\em E-mail addresses:} daniel@sunam1.mat.ucm.es,
azagra@ccr.jussieu.fr, \\ marjim@sunam1.mat.ucm.es,
marjim@ccr.jussieu.fr

\end{document}